\documentclass[preprint,3p]{elsarticle}

\journal{Elsevier}
\biboptions{sort&compress}

\usepackage[utf8]{inputenc}
\usepackage{bm}
\usepackage{amsmath}
\usepackage[hidelinks]{hyperref}
\hypersetup{colorlinks=true,linkcolor=blue,urlcolor=blue}
\usepackage{multicol}
\usepackage{booktabs}
\usepackage{amsthm}
\usepackage[table,xcdraw,dvipsnames]{xcolor}
\usepackage{graphicx}
\usepackage{float}
\usepackage{subcaption}
\usepackage{algorithm}
\usepackage{algpseudocode}
\usepackage{pdflscape}
\usepackage{array}
\usepackage{multirow}
\usepackage{diagbox}
\usepackage{cleveref}
\usepackage{amsfonts}
\usepackage{amssymb, comment}
\usepackage{ulem}
\usepackage{fourier-orns}
\usepackage[labelfont=bf,labelsep=period]{caption}
\usepackage[labelfont=rm,font=footnotesize]{subcaption}

\numberwithin{equation}{section}

\DeclareMathAlphabet{\altmathcal}{OMS}{cmsy}{m}{n}
\DeclareMathAlphabet{\altmathcalb}{OMS}{cmsy}{b}{n}
\DeclareMathAlphabet{\mathcalboondox}{U}{BOONDOX-calo}{m}{n}
\DeclareMathAlphabet{\mathbbmsl}{U}{bbm}{m}{sl}

\newdefinition{remark}{Remark}

\newcommand{\orcid}[1]{\href{https://orcid.org/#1}{\texorpdfstring{\includegraphics*[width=8pt]{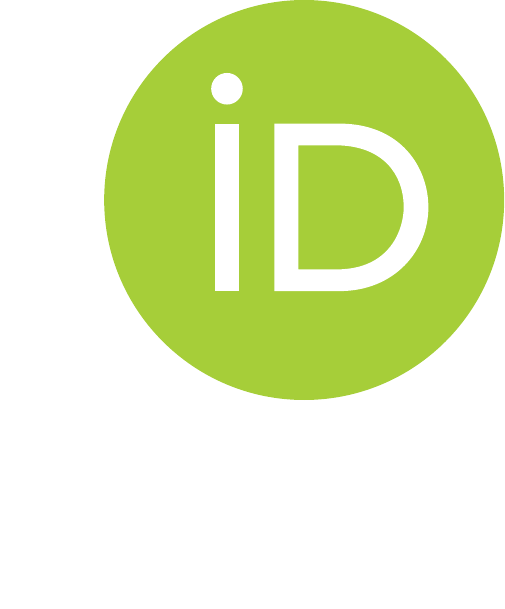}}~~}{}}

\definecolor{drot}{rgb}{0.7,0,0.1}

\begin{document}
\baselineskip14pt

\begin{frontmatter}

\title{Efficient Numerical Integration for Finite Element Trunk Spaces in 2D and 3D using Machine Learning: A new Optimisation Paradigm to Construct Application-Specific Quadrature Rules}

\author[bcam]{Tomas Teijeiro\orcid{0000-0002-2175-7382}\corref{cor1}}
\cortext[cor1]{Corresponding author}
\ead{tteijeiro@bcamath.org}
\author[affil]{Pouria Behnoudfar\orcid{0000-0003-4301-3728}}
\author[cunef]{Jamie M. Taylor\orcid{0000-0002-5423-828X}}
\author[upv,bcam,ikerbasque]{David Pardo\orcid{0000-0002-1101-2248}}
\author[curtin]{Victor M. Calo\orcid{0000-0002-1805-4045}}

\address[bcam]{{BCAM -- Basque Center for Applied Mathematics}, {Bilbao, Basque Country, Spain}}

\address[affil]{{University of Wisconsin-Madison, Department of Mathematics, Madison, WI, USA}}

\address[cunef]{{Department of Mathematics, CUNEF Universidad}, {Madrid, Spain}}
\address[upv]{{University of the Basque Country (UPV/EHU)}, {Leioa, Basque Country, Spain}}
\address[ikerbasque]{{Ikerbasque -- Basque Foundation for Sciences}, {Bilbao, Basque Country, Spain}}
\address[curtin]{{School of Electrical Engineering, Computing and Mathematical Sciences, Curtin University}, {Perth, Australia}}

\begin{abstract}

Finite element methods usually construct basis functions and quadrature rules for multidimensional domains via tensor products of one-dimensional counterparts. While straightforward, this approach results in integration spaces larger than necessary, especially as the polynomial degree \( p \) or the spatial dimension increases, leading to considerable and avoidable computational overhead. This work starts from the hypothesis that reducing the dimensionality of the polynomial space can lead to quadrature rules with fewer points and lower computational cost, while preserving the exactness of numerical integration. We use \textit{trunk spaces} that exclude high-degree monomials that do not improve the approximation quality of the discrete space. These reduced spaces retain sufficient expressive power and allow us to construct smaller (more economical) integration domains. Given a maximum degree \( p \), we define trial and test spaces \( U \) and \( V \) as 2D or 3D trunk spaces and form the integration space \( \mathcal{S} = U \otimes V \). We then construct exact quadrature rules by solving a non-convex optimisation problem over the number of points \( q \), their coordinates, and weights. We use a shallow neural network with linear activations to parametrise the rule, and a random restart strategy to mitigate convergence to poor local minima. When necessary, we dynamically increase \( q \) to achieve exact integration. Our construction reaches machine-precision accuracy (errors below \( 10^{-22} \)) using significantly fewer points than standard tensor-product Gaussian quadrature: up to 30\% reduction in 2D for \( p \leq 10 \), and 50\% in 3D for \( p \leq 6 \). These results show that combining the mathematical understanding of polynomial structure with numerical optimisation can lead to a practical and extensible methodology for improving the adaptiveness, efficiency, and scalability of quadrature rules for high-order finite element simulations.

\end{abstract}


\end{frontmatter}



\section{Introduction}
\label{sec:Introduction}

Partial differential equations (PDEs) constitute the mathematical foundation for describing various physical processes across science and engineering. They model the behaviour of phenomena such as heat transfer, fluid flow, electromagnetic fields, and structural deformation. Due to the complexity of real-world systems, analytical solutions to these equations are seldom available, making numerical approximation methods indispensable. Among these, the Finite Element Method (FEM) is a widely adopted framework for numerically solving PDEs~\cite{Ern:2004, Hughes:2000, Demkowicz:2006, Pardo:2021Book, Liu2022}. FEM decomposes the domain into small subdomains (finite elements) over which the solution is approximated using locally defined basis functions, typically polynomials. This method offers several advantages: it is backed by a rigorous theoretical framework, it accommodates diverse types of equations, it adapts flexibly to complex geometries, and it scales efficiently with high-performance computing architectures.  In a forward problem, the objective is to determine the system's state, given the initial and boundary conditions, governing equations, and physical parameters~\cite{Huebner2001, Demkowicz:2007, Cottrell:2009, Pardo:2021Book}. The discrete approximations generated by FEM enable precise simulation, leading to predictive modelling capabilities and the optimisation of engineering designs. 

Nevertheless, one of the primary limitations of FEM lies in its high computational demands~\cite{Liu2022}. Numerous strategies have been proposed to address this issue. Among them are high-order methods and adaptive mesh refinement. High-order FEM, in particular, delivers exponential convergence with respect to the polynomial degree \( p \), assuming the underlying solution is sufficiently smooth. However, this increased accuracy is accompanied by polynomial growth in computational cost, which can become prohibitive. In \( p \)-FEM, the polynomial degree can vary from one element to another and along different coordinate directions~\cite{ Babuska:1981}. Such anisotropic variation introduces difficulties in maintaining global continuity across element interfaces. When adjacent elements possess different polynomial degrees in the tangential direction, only the lower-degree basis functions can be matched across the interface to ensure global \( C^0 \)-conformity. The combined \( hp \)-FEM strategy has been extensively studied and extended to three-dimensional settings. It has proven effective in the presence of geometric complexity, singularities, and vector-valued equations such as Maxwell's equations~\cite{Demkowicz:1989, Demkowicz:2006, Demkowicz:2007}. This approach exploits local adaptivity and exponential convergence, supported by \textit{a posteriori} error estimation and residual minimisation techniques. In regions where the solution exhibits sharp gradients or interfaces, the mesh is algorithmically refined, while coarser resolution is retained elsewhere~\cite{ Demkowicz:2010, Carstensen:2016, Cier:2021, Giraldo:2023, Hasbani:2023}; these techniques are robust in various challenging applications~\cite{ Poulet:2023, Labanda:2022, Demkowicz:2024, Henneking:2023, Henneking:2021}.

An additional line of improvement has been the development of multi-level \( hp \)-methods~\cite{DiStolfo:2016, Zander:2016, Zander:2015, Zander:2022}. These methods avoid complex hierarchical data structures by decoupling coarse and refined representations. Polynomial bases are constructed independently at each level, simplifying the assembly process and enhancing computational scalability. The resulting frameworks exhibit efficient performance during the assembly of linear systems and are freely available as software libraries~\cite{Kopp:2022a, Jomo:2019, Jomo:2021}. These tools have been successfully applied in diverse fields, including biomedicine~\cite{Elhaddad:2018} and additive manufacturing~\cite{Kopp:2022b}.

Trunk spaces~\cite{ Szabo2011} reduce the computational cost by reducing the number of functions used to approximate the solution discretely while preserving the approximation quality of the space. These are reduced polynomial spaces derived by omitting higher-order monomials from the full tensor-product basis. Specifically, instead of employing all basis functions in the full multidimensional space, one retains only those satisfying
\[
\sum_{i=0}^{d-1} r_i \leq \max(p_i),
\]
where \( r_i \) is the polynomial degree, and \( p_i \) is the maximum degree permitted in the \( i \)-th coordinate direction. This restriction significantly reduces the number of degrees of freedom, while preserving the convergence rate of the finite element discretisation. Although monomials themselves are not appropriate for constructing globally compatible shape functions, we can translate the reduced basis concept to compatible polynomial families, such as integrated Legendre polynomials~\cite{Demkowicz:2006, Demkowicz:2007} or classical Lagrange polynomials~\cite{ Hughes:2000}, which retain interelement continuity and conformity while benefiting from the reduced dimensionality of the trunk space.

The reduction of degrees of freedom entails a trade-off. For a fixed polynomial order, trunk spaces may exhibit a minor loss in pointwise accuracy. However, this can often be compensated by increasing the polynomial degree. Reducing the number of active basis functions leads to significant memory and computational time savings. The remaining degrees of freedom exhibit broader coupling, potentially increasing the matrix assembly and solution cost. As dimensionality increases, the advantages of trunk spaces become even more pronounced. In \( d \) dimensions, the number of unknowns per element in a trunk space is approximately \( d! \) times smaller than in the corresponding tensor-product space. Geometrically, this scaling arises from comparing the volume of a unit simplex, which is \( 1/d! \), to that of a unit hypercube, which is one. This contrast reveals the efficiency of using a simplex-like structure in high dimensions. These motivations are in line with contemporary data-driven and machine-learning approaches that aim to reduce numerical burden. Examples include deep learning~\cite{Hashash2004}, physics-informed neural networks~\cite{Karniadakis2021}, reduced-order models~\cite{Calo:2016, Lu2021}, hierarchical neural networks~\cite{Saha2021}, and hybrid strategies combining several of these methods~\cite{Efendiev:2009, Ghommem:2013, Ghommem:2014, Zhang2022}. However, all such frameworks depend critically on the efficient numerical integration of the discretised system.

Herein, we address this foundational operation (numerical integration) by constructing specialised quadrature rules tailored to trunk spaces. Unlike conventional Gaussian quadrature built from tensor-product rules, we target integration rules that are exact on reduced spaces and require fewer evaluation points. Since no analytical method exists to identify optimal rules of minimal cardinality, we frame this as a non-convex optimisation problem. Using a loss function to measure integration error, we use machine learning techniques to search for rules that meet the desired properties. The optimisation uses shallow neural networks with random restarts and adaptive strategies to increment the number of quadrature points when necessary.

The main contributions of this study are:
\begin{itemize}
    \item A collection of precomputed quadrature rules for efficient and exact integration on 2D and 3D trunk product spaces with varying polynomial degrees.  
    \item  A general methodology to find new quadrature rules tailored to custom spaces and applications, such as those with different polynomial degrees in each dimension.
\end{itemize} 

In high-order finite element methods, integration cost escalates rapidly with problem complexity, particularly in nonlinear and high-dimensional contexts. Evaluating shape functions and their derivatives across numerous integration points significantly contributes to computational cost. While tensor-product rules are optimal for full polynomial spaces, they are unnecessarily expensive when applied to trunk spaces~\cite{Yosibash1994}. High-order tensor-product bases introduce spectral artefacts such as stopping bands and optical branches~\cite{Hughes:2008, Barton:2018, Puzyrev:2017}, further motivating reduced space use. Simplifying the integration domain enables the construction of quadrature rules with fewer points while preserving exactness. This reduces computational expense without compromising accuracy.

The structure of this manuscript is as follows: 
\autoref{sec:ProblemFormulation} introduces the mathematical formulation of the integration problem, including the definition of function spaces and the loss function. 
\autoref{sec:MLAlgorithm} presents the optimisation algorithm, including the use of gradient descent with random restarts. 
\autoref{sec:Results} reports numerical experiments for 2D and 3D trunk spaces, evaluating the integration error and comparing point counts with standard Gaussian quadrature. 
\autoref{sec:CaseStudies} illustrates two practical applications for employing these rules. 
Finally, \autoref{sec:Conclusions} summarises our findings and outlines future research directions.


\section{Exact Quadrature Rules for Trunk Spaces}
\label{sec:ProblemFormulation}

\subsection{Trunk Spaces}
\label{sec:trunk}

We develop exact quadrature rules for products of trunk spaces---also known as serendipity spaces---that require fewer evaluations than standard methods based on tensor-product Gaussian quadrature. Trunk spaces consist of all monomials whose \textit{total degree} does not exceed a specified maximum, in contrast to tensor-product spaces, which include all combinations of degrees in each coordinate direction up to the maximum per dimension~\cite[Section 5.2.1]{ Szabo2011}.

While tensor-product spaces are simple to construct and analyse, they include many high-order terms that contribute minimally to the solution's precision but significantly increase the number of required quadrature points, and also the size (and cost) of corresponding linear systems. By constraining our function space to trunk spaces, we can formulate more efficient quadrature rules—especially in two- and three-dimensional cases involving quadrilateral or hexahedral finite elements. These rules reduce the numerical cost, particularly during the assembly of stiffness and mass matrices, while maintaining convergence rates comparable to full tensor-product spaces.

Since trunk spaces are strict subsets of their corresponding tensor-product spaces, they offer a more compact basis for high-order finite element integration. Although trunk spaces can reduce pointwise accuracy for a fixed polynomial degree and mesh, we can recover the accuracy by increasing the polynomial degree. This typically results in fewer degrees of freedom overall, although each unknown has higher connectivity and influences more neighboring elements~\cite{ Kopp:2022a}. As explained in~\autoref{sec:Introduction}, the efficiency gains trunk spaces offer become increasingly pronounced as the number of spatial dimensions increases.

\subsubsection{Illustrative Example in 2D}

Figure~\ref{fig:trunk} illustrates a 2D trunk space for a quadrilateral element with polynomial degree $p = 3$. We display the monomials of the full tensor-product space (a product of two 1D degree-$p$ polynomial spaces) and indicate the monomials removed to form the trunk space, i.e., those whose total degree exceeds $p$. Two higher-order functions, $x^p y$ and $x y^p$, are retained despite their total degree being $p+1$, to preserve approximation quality at interfaces and ensure smooth variation into the element interior~\cite[Section 5.2.1]{Szabo2011}.

\begin{gather*}
    1\\
    x\quad~y\\
    x^2\quad~xy\quad~y^2\\
    x^3\quad~x^2y\quad~xy^2\quad~y^3\\
    \hbox{$x^3y$\quad~\sout{$x^2y^2$}\quad~$xy^3$}\\
    \hbox{\sout{$x^3y^2$}\quad~\sout{$x^2y^3$}}\\
    \hbox{\sout{$x^3y^3$}}
\end{gather*}
\captionof{figure}{\label{fig:trunk}Monomials in a 2D trunk space of degree $p = 3$ on a quadrilateral element. Removed terms (those of excessive total degree) are shown crossed out.}

\subsubsection{Target Integration Space}

We define the space on which the quadrature rule must be exact. For a given degree \( p \), we define both the trial space \( U \) and the test space \( V \) as trunk spaces of degree \( p \) on the reference domain \( [0,1]^2 \) in 2D (and \( [0,1]^3 \) in 3D). The set of functions requiring exact integration lies in the span of the product space
\[
\mathcal{S} = \text{span}(U \otimes V).
\]
This space contains all bilinear combinations of basis functions from the trial and test spaces. Also, as the gradients of elements of $\mathcal{S}$ are within $\mathcal{S}$ as well, exact integration of the mass matrix implies exact integration of the stiffness matrix. The cardinality of this space, denoted \( |\mathcal{S}|(p) \), is essential for determining the number of quadrature points required to guarantee exact integration. The number of basis functions in \( \mathcal{S} \) determines the number of degrees of freedom in the corresponding discrete bilinear form. To calculate $|\mathcal{S}|(p)$, we take as a preliminary the number of basis functions in trunk spaces of degree $p$, which we denote $|\mathcal{T}|(p)$. It is well known~\cite{Duster2001, Szabo2011} that $|\mathcal{T}|(p)$ follows these equations for the 2D and 3D cases:

\begin{alignat}{2}
|\mathcal{T}|_{2D}(p) &=  \frac{(p+1)(p+2)}{2} + 2, &\quad \text{for } p \geq 2 \label{eq:trunk_size_2D}\\
|\mathcal{T}|_{3D}(p) &= \frac{(p+1)(p+2)(p+3)}{6} + 3p + 3, &\quad \text{for } p \geq 3 \label{eq:trunk_size_3D}
\end{alignat}

We also establish lower and upper bounds for the product space \( \mathcal{S} \):
\begin{alignat}{2}
|\mathcal{T}|_{2D}(2p) &~<~ |\mathcal{S}|_{2D}(p) &~<~& (2p + 1)^2 \label{eq:size_limits_2D}\\
|\mathcal{T}|_{3D}(2p) &~<~ |\mathcal{S}|_{3D}(p) &~<~& (2p + 1)^3 \label{eq:size_limits_3D}
\end{alignat}

These inequalities express that the span of the product of two trunk spaces is strictly larger than a trunk space of degree \( 2p \), but still smaller than the full tensor-product space of that degree. This property is evident from the earlier example in Figure~\ref{fig:trunk}. For instance, in two dimensions with \( p = 3 \), the trunk space \( \mathcal{T}_{2D}(6) \) contains monomials of total degree up to 7. The product space \( \mathcal{S}_{2D}(3) \), formed from degree-3 trunk functions, contains monomials up to degree 8. Meanwhile, the full tensor-product space at degree 3 spans monomials up to degree 12. From these bounds, we infer that \( |\mathcal{S}|_{2D}(p) \) scales quadratically in \( p \), while \( |\mathcal{S}|_{3D}(p) \) exhibits cubic growth. We compute explicit expressions for these dimensions by constructing product spaces for representative values of \( p \), and applying polynomial regression to obtain:
\begin{alignat}{2}
|\mathcal{S}|_{2D}(p) &=  2p^2 + 5p + 4, &\quad \text{for } p \geq 2 \label{eq:size_2D}\\
|\mathcal{S}|_{3D}(p) &= \frac{4 p^3 + 24 p^2 + 56 p + 21}{3}, &\quad \text{for } p \geq 3 \label{eq:size_3D}
\end{alignat}
These estimates guide our assumptions about the minimum number of quadrature points needed to integrate all elements of \( \mathcal{S} \) exactly, and they are critical for our optimisation framework.

\section{Finding Quadrature Rules as an Optimisation Problem}
\label{sec:optimization_problem}

This section constructs quadrature rules for integrating functions in the space $\mathcal{S}$ using the smallest possible number of points. Instead of relying on classical quadrature techniques, we approach the task as a numerical optimisation problem. A quadrature rule is a collection of weighted points approximating definite integrals. In two dimensions, we define the quadrature rule as a set of $q$ triplets $(x_i, y_i, w_i)$, where $(x_i, y_i)$ are the coordinates of the $i$-th integration point, and $w_i$ is its associated weight. For three-dimensional problems, the set expands to include the $z$ coordinate, yielding quadruplets $(x_i, y_i, z_i, w_i)$. Explicitly, we write:
\begin{equation}
Q := 
\begin{cases}
\displaystyle
 \{(x_i, y_i, w_i)\}_{i=1}^q, & \text{in 2D},\\[0.5em]
\displaystyle
\{(x_i, y_i, z_i, w_i)\}_{i=1}^q, & \text{in 3D}.
\end{cases}
\end{equation}
such that the quadrature integrates all functions in $\mathcal{S}$ exactly $\forall f \in \mathcal{S}$:
\begin{equation}
\label{eq:quadrature}
\begin{aligned}
\iint_{[0,1]^2} f(x, y)\, dx\,dy &= \sum_{i=1}^q f(x_i, y_i)\, w_i \qquad &\text{(2D case)},\\[0.5em]
\iiint_{[0,1]^3} f(x, y, z)\, dx\,dy\,dz &= \sum_{i=1}^q f(x_i, y_i, z_i)\, w_i \qquad &\text{(3D case)}.
\end{aligned}
\end{equation}

We further require that all coordinates and weights fall strictly within the unit domain:
\[
\begin{cases}
x_i, y_i, w_i \in (0, 1), & \text{for 2D},\\[0.5em]
x_i, y_i, z_i, w_i \in (0, 1), & \text{for 3D},
\end{cases}
\quad \text{and} \quad
\sum_{i=1}^q w_i = 1.
\]
We estimate the smallest viable value for $q$ by observing that each function $f$ in $\mathcal{S}$ yields one equation when imposed in equation~\eqref{eq:quadrature}. The dimension of $\mathcal{S}$, denoted $|\mathcal{S}|(p)$, equals the number of such constraints. Meanwhile, each integration point contributes several degrees of freedom: three in 2D ($x$, $y$, $w$) and four in 3D ($x$, $y$, $z$, $w$). Thus, a dimension-counting argument suggests that the lower bound for the number of required quadrature points is:
\[
q = \left\lceil \frac{|\mathcal{S}|_{2D}(p)}{3} \right\rceil, \qquad
q = \left\lceil \frac{|\mathcal{S}|_{3D}(p)}{4} \right\rceil.
\]

\subsection{Comparison with Tensor-Product Gaussian Rules}

Standard Gaussian quadrature requires $q' = (p+1)^2$ points in 2D and $q' = (p+1)^3$ in 3D. Our method can potentially reduce this number by:

\begin{equation}
\text{Savings} = 1 - \frac{q}{q'}.
\label{eq:savings}
\end{equation}

By inserting the closed-form expressions for $|\mathcal{S}|(p)$ provided in equations~\eqref{eq:size_2D} and~\eqref{eq:size_3D}, and examining the asymptotic behaviour as $p \to \infty$, we observe that the ideal savings converge to $\frac{1}{3}$ in 2D and $\frac{2}{3}$ in 3D.

\subsection{Loss Function for Optimisation}

We frame the optimisation problem precisely. Let $\{v_i(x, y)\}_{i=1}^n$ be a basis for $\mathcal{S}$, where $n = |\mathcal{S}|$. We define the Gram matrix $M$ with entries:
\[
M_{ij} := 
\begin{cases}
\displaystyle
\iint_{[0,1]^2} v_i(x, y) v_j(x, y) \, dx \, dy, & \text{in 2D},\\[1.5em]
\displaystyle
\iiint_{[0,1]^3} v_i(x, y, z) v_j(x, y, z) \, dx \, dy \, dz, & \text{in 3D}.
\end{cases}
\]
and the integration error for each basis function:
\[
e(v_i) := 
\begin{cases}
\displaystyle
\left| \iint_{[0,1]^2} v_i(x, y)\, dx\,dy - \sum_{j=1}^q v_i(x_j, y_j)\, w_j \right|, & \text{2D},\\[0.5em]
\displaystyle
\left| \iiint_{[0,1]^3} v_i(x, y, z)\, dx\,dy\,dz - \sum_{j=1}^q v_i(x_j, y_j, z_j)\, w_j \right|, & \text{3D}.
\end{cases}
\]

We seek to minimise the worst-case relative integration error across all functions in $\mathcal{S}$. The loss function used in the optimisation is:
\begin{equation}
\label{eq:loss}
L(Q) := \max_{f \in \mathcal{S}} \frac{|e(f)|}{\|f\|_{L^2}} 
= \left( \sum_{i,j=1}^n (M^{-1})_{ij} \, e(v_i) \, e(v_j) \right)^{1/2}.
\end{equation}

A quadrature rule $Q$ is exact on $\mathcal{S}$ if and only if $L(Q) = 0$. Apart from the dimensionality of the integration points and the basis size, the optimisation problem formulation remains identical for 2D and 3D geometries.

Our formulation exploits several methodological and computational advantages of using trunk spaces in high-order finite element methods, including the preservation of the asymptotic convergence rate and the factorial reduction of the cardinality with respect to dimension.

On the other hand, the target space for integration is not a trunk space, although it is constructed as a product of trunk spaces. Also, the implicit interplay between the number of basis functions and the minimum number of quadrature points requires accurate and robust estimates. Each quadrature point contributes multiple degrees of freedom (three in 2D, four in 3D), and the minimal count of such points is bounded from below by the ratio of the dimension of the target space to this contribution. Crucially, this relation provides a theoretical baseline for optimality, enabling direct comparison to Gaussian quadrature and thus giving an upper bound for achievable efficiency gains. Finally, the optimisation landscape is highly non-convex, which complicates the minimisation process and requires robust strategies to avoid local minima. The solution adopts a machine-learning-inspired algorithm incorporating adaptive gradient descent with restart and early-stopping mechanisms. This setup allows the method to iteratively improve candidate quadrature rules until the loss function, which quantifies the integration error across the basis, reaches a numerically negligible value. The algorithm further allows for adaptivity in the number of quadrature points, relaxing the minimality constraint if an exact rule proves elusive. These procedural choices are not merely technical artefacts but essential responses to the intrinsic complexity of the quadrature discovery problem in trunk-product spaces.


\section{Machine Learning-Based Optimisation Strategy}
\label{sec:MLAlgorithm}

The core challenge of solving the optimisation problem defined in Section~\ref{sec:ProblemFormulation} lies in the highly \textit{non-convex} nature of the loss function, see~\eqref{eq:loss}. As a result, standard gradient descent methods may frequently become trapped in local minima, failing to find an exact quadrature rule. Nevertheless, the specific structure of our problem offers a distinct advantage. Provided that the number of quadrature points $q$ exceeds a theoretical minimum (i.e., $q > \frac{|\mathcal{S}|_{2D}}{3}$ in 2D, or $q > \frac{|\mathcal{S}|_{3D}}{4}$ in 3D), we observed that there exists not just a unique solution but rather an \textit{infinite multiplicity} of exact quadrature rules. This heuristic property increases the likelihood of success in finding a global minimum through numerical optimisation.

Moreover, in contrast to many traditional machine learning applications where ground truth is inaccessible or ill-defined, our algorithm can immediately verify whether a candidate quadrature rule is exact by evaluating whether the loss is below a specified numerical tolerance. This property simplifies convergence assessment and transforms the problem into a well-posed verification task once a suitable candidate has emerged. These observations allow us to employ a robust optimisation approach inspired by machine learning: \textbf{restarted gradient descent} in conjunction with \textbf{early stopping}; a combination that balances computational efficiency and robustness against entrapment in local minima.

\subsection{Algorithm Overview}

\begin{algorithm}[h!]
\caption{Restarted gradient descent for discovering exact quadrature rules in trunk-product approximation spaces}
\label{alg:restarted_gdescent}

\algrenewcommand{\alglinenumber}[1]{#1}
\small
\begin{algorithmic}[1]
    \State \textbf{Input}: Desired polynomial degree $p$, spatial dimension $d \in \{2,3\}$
    \State \textbf{Output}: Quadrature rule $Q$ such that $L(Q) < 10^{-22}$
        \vspace{0.5em}
        
    \State \textbf{Initial Setup:}
        \State Define the optimiser $\gamma$ as the Yogi algorithm with learning rate set to $10^{-2}$ \label{line:optimizer}
        
        \If{$d = 2$}
            \State Construct trunk spaces $U = V$ of degree $p$ defined over the domain $[0,1]^2$
            \State Define the integration space $\mathcal{S} = \text{span}(U \otimes V)$, with basis functions $\{v_i(x,y)\}_{i=1}^n$
            \State Estimate the minimal number of quadrature points $q \gets \lceil \frac{n}{3} \rceil$
        \ElsIf{$d = 3$}
            \State Construct trunk spaces $U = V$ of degree $p$ on $[0,1]^3$
            \State Define $\mathcal{S} = \text{span}(U \otimes V)$, with basis $\{v_i(x,y,z)\}_{i=1}^n$
            \State Estimate $q \gets \lceil \frac{n}{4} \rceil$
        \EndIf \label{line:q}
        
    \State Precompute: the inverse Gram matrix $M^{-1}$ and the exact integral of each $v_i$ over $[0,1]^d$

        \Statex

    \State \textbf{Outer Restart Loop:} Attempt multiple optimisation runs, increasing redundancy if needed
    \Loop \label{line:outer}
        \For{$iter = 1$ to $10^4$} \label{line:iters}

            \If{$d = 2$}
                \State Initialise quadrature rule $Q_0 = \{(x_i, y_i, w_i)\}_{i=1}^q$ 
                \State Each parameter $x_i, y_i, w_i$ sampled uniformly: $x_i, y_i, w_i \sim \mathcal{U}(0,1)$ \label{line:rand2}
            \ElsIf{$d = 3$}
                \State Initialise $Q_0 = \{(x_i, y_i, z_i, w_i)\}_{i=1}^q$
                \State Each parameter sampled: $x_i, y_i, z_i, w_i \sim \mathcal{U}(0,1)$ \label{line:rand3}
            \EndIf

            \State Normalise the weights so they sum to unity: $w_i \gets w_i / \sum_j w_j$ \label{line:norm}
                
            \Statex

            \State \textbf{Inner Optimisation Loop:} Improve quadrature rule using gradient descent
            \For{$e = 1$ to $10^6$} \label{line:inner}
                \State Update rule using optimiser: $Q_e \gets \gamma(Q_{e-1}, \nabla L(Q_{e-1}))$ \label{line:gdescent}

                \If{the loss satisfies $L(Q_e) < 10^{-22}$} \label{line:exact_rule}
                    \State \Return $Q_e$ \Comment The quadrature rule is considered exact
                \ElsIf{no significant improvement over last 100 steps:}
                    \State EarlyStopping criterion triggered: 
                    \Statex \hspace{2em} $\text{EarlyStopping}(L(Q_{e-100}), \ldots, L(Q_e), 10^{-23})$ \label{line:early_stop}
                    \State \textbf{break} \Comment Convergence stagnated
                \EndIf
            \EndFor \label{line:end_inner}

        \EndFor

        \State Increment the number of quadrature points: $q \gets q + 1$ 
        \Comment Expands the search space by relaxing the constraint
    \EndLoop
\end{algorithmic}
\end{algorithm}

Algorithm~\ref{alg:restarted_gdescent} encapsulates the detailed procedure. The method takes the polynomial degree $p$ as input, which governs the definition of the trunk spaces $U$ and $V$. The algorithm outputs the first found quadrature rule $Q$ that satisfies the accuracy condition by bringing the loss below the threshold.

\paragraph{Initialisation} Lines~\ref{line:optimizer}--\ref{line:q} perform preliminary setup. We construct the trunk spaces $U = V$ of degree $p$ and form the integration space $\mathcal{S} = \text{span}(U \otimes V)$. Then, we compute the inverse of the Gram matrix $M^{-1}$ associated with the basis of $\mathcal{S}$, alongside the vector of exact integrals for each basis function $v_i \in \mathcal{S}$. We estimate the theoretical lower bound for the number of integration points $q$ as:
    \[
    q = 
    \begin{cases}
    \left\lceil \frac{|\mathcal{S}|_{2D}(p)}{3} \right\rceil, & \text{in two dimensions,}
    \\[1.5em]
    \left\lceil \frac{|\mathcal{S}|_{3D}(p)}{4} \right\rceil, & \text{in three dimensions.}
    \end{cases}
    \]

\paragraph{Outer Optimisation Loop (Restart Strategy)} The algorithm runs an outer loop, see line~\ref{line:outer}, which systematically performs multiple optimisation attempts. In each attempt, or \textit{restart}, we generate a new quadrature rule $Q_0$  by sampling $q$ spatial coordinates (either $(x_i, y_i)$ in 2D (line~\ref{line:rand2}) or $(x_i, y_i, z_i)$ in 3D (line~\ref{line:rand3})) from uniform probability distributions. Initial weights $w_i$ are normalised to satisfy a partition-of-unity condition (line~\ref{line:norm}).

\paragraph{Inner Loop (Gradient-Based Optimisation)} We iteratively refine the quadrature configuration using a gradient-based method within each restart. Specifically, we employ the Yogi optimiser~\cite{NEURIPS18_adaptive}, an adaptive gradient descent algorithm that dynamically adjusts learning rates based on recent gradient statistics, promoting stability and convergence even under non-convex conditions. The learning rate, a key hyperparameter, is fixed empirically at $10^{-2}$ (line~\ref{line:optimizer}).

\paragraph{Termination Criteria:} The inner loop concludes under any of the following scenarios:
\begin{enumerate}
    \item A quadrature rule is deemed exact if the loss $L(Q)$ drops below $10^{-22}$ (line~\ref{line:exact_rule}). Given numerical precision constraints, this value is a practical surrogate for zero. \label{en:prec}
    \item Reaching the maximum iteration number (set to $10^6$) without reaching the zero surrogate (line~\ref{line:inner}).
    \item Optimiser stagnates, observing negligible loss reduction (less than $10^{-23}$) across the last 100 iterations. This constitutes an \textit{early stopping} trigger (line~\ref{line:early_stop}), a common strategy in model training~\cite{Raschka2020}.
\end{enumerate}
    
\paragraph{Adaptivity via Restarting} Should the algorithm fail to yield a valid quadrature rule after $10^4$ restarts, the value of $q$ is incremented by one (line~\ref{line:outer}). This relaxes the problem by further underdetermining the system, thereby increasing the feasibility of locating an exact solution in subsequent trials.

This combination of stochastic initialisation, adaptive descent, and progressive relaxation provides a flexible yet effective pathway toward discovering exact integration formulas over trunk-product spaces, even in severe non-convexities and high-dimensional parameter domains.


\section{Numerical Results}
\label{sec:Results}

\subsection{Implementation Details}

We use a \texttt{Python} implementation of the optimisation framework, with heavy reliance on the \texttt{JAX} library~\cite{ jax2018github} for defining and evaluating basis functions, as well as for computing gradients of the loss function through automatic differentiation. The Yogi optimiser drives the optimisation loop, as implemented in the \texttt{Optax} library~\cite{ deepmind2020jax}, a well-tested tool in modern machine learning applications.

We run all the numerical experiments on a standard laptop PC running Linux kernel version 6.8, with an Intel(R) Core(TM) i7-1250U CPU and 32~GB of RAM. The Python environment used for testing was version 3.12.

Although Section~\ref{sec:trunk} introduces trunk spaces based on monomial sets, we did not apply them directly due to their intrinsic numerical instability. Monomial bases, while conceptually transparent, tend to produce Gram matrices with high condition numbers, compromising the accuracy and efficiency of subsequent computations. Instead, we use a more robust hierarchical basis constructed from tensor products of Legendre-type shape functions, following the methodology outlined in~\cite{ Duster2001}. This choice improves numerical stability without sacrificing representational richness.

We make the source code publicly available to promote reproducibility and encourage future development. In addition to the implementation, we provide exact quadrature rules in double precision, covering 2D trunk spaces up to $p = 10$ and 3D trunk spaces up to $p = 6$. These resources and supplementary documentation are freely accessible at: \url{https://github.com/Mathmode/trunk-spaces-quadrature}.

\subsection{General Experiments}
\label{sub:GeneralExperiments}

We evaluate the effectiveness and scalability of the proposed method by conducting extensive experiments targeting exact quadrature rules for polynomial spaces defined over single quadrilateral (2D) and hexahedral~(3D) elements. The objective for each polynomial degree $p$ was to identify at least one quadrature rule that integrates the space exactly. Even though the space of exact rules is generally infinite (as Section~\ref{sec:MLAlgorithm} discusses), only a single valid rule is required per configuration.

\begin{figure}[t]
     \centering
     \begin{subfigure}[b]{0.48\textwidth}
         \centering
         \includegraphics[width=\textwidth]{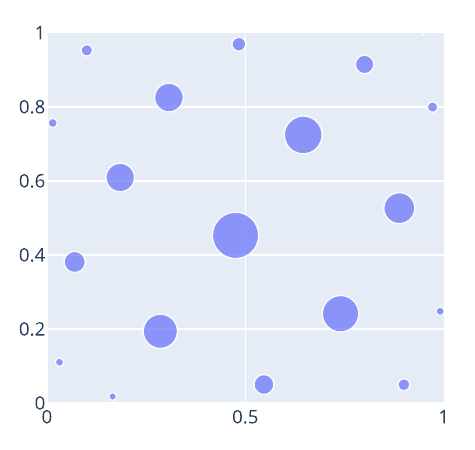}         
     \end{subfigure}
     \hfill
     \begin{subfigure}[b]{0.48\textwidth}
         \centering
         \includegraphics[width=\textwidth]{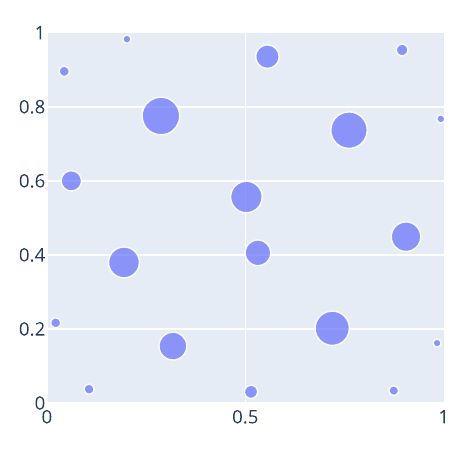}
     \end{subfigure}
        \caption{Two distinct quadrature rules for a 2D trunk space with $p = 4$. Each circle represents a quadrature node at $(x_i, y_i)$, with area scaled to its weight $w_i$. Both rules achieve exact integration.}
        \label{fig:rule_examples}
\end{figure}

Figure~\ref{fig:rule_examples} illustrates two distinct yet valid exact quadrature rules for the case $p = 4$ in two dimensions. This visual comparison underscores the inherent non-uniqueness of the solution space: multiple point configurations, each with different node locations and weight distributions and unrelated by symmetry, may satisfy the same exactness condition. Each quadrature point is shown as a circle whose position corresponds to the integration node $(x_i, y_i)$, and whose area is proportional to the weight $w_i$.

Table~\ref{tab:IntegrationSavings} summarizes our numerical findings. For each polynomial degree and spatial dimension, the table compares:
\begin{itemize}
    \item the number of quadrature points required by classical Gaussian tensor-product rule ($q'$),
    \item the smallest number of points in our discovered exact rule ($q$),
    \item the relative reduction in quadrature size, expressed as a percentage,
    \item and the number of optimisation restarts necessary to obtain the result (as per line~\ref{line:iters} in Algorithm~\ref{alg:restarted_gdescent}).
\end{itemize}

In 2D, results demonstrate that exact quadrature rules can be reliably obtained for polynomial degrees up to $p=10$, with substantial savings in the number of integration points, approaching the theoretical lower bounds for $p \geq 3$. One exception is the case $p = 9$ requiring a large number of restarts (3072), likely due to unfavorable random initialization rather than algorithmic failure.

In 3D, the situation becomes significantly more complex. While rules were successfully identified up to $p = 6$, the number of required restarts increased sharply with polynomial degree. Moreover, for $p = 5$ and $p = 6$, the discovered rules used slightly more points than the known theoretical optima (117 and 174, respectively). Nonetheless, the computational gains are evident: for these cases, our quadrature rules employ nearly 50\% fewer points than Gaussian tensor-product. This efficiency level has meaningful implications for large-scale simulations in finite element analysis, where quadrature cost is often a dominant factor.

\begin{table}[!h]\centering
    \caption{Summary of optimised quadrature results compared to standard Gaussian tensor-product rules. Values shown include polynomial degree $p$, the number of integration points in Gaussian rules ($q'$), our optimised rules ($q$), relative savings, and the number of optimisation restarts needed. All data corresponds to a single quadrilateral (2D) or hexahedral (3D) element.}
	\label{tab:IntegrationSavings}
    \small
	\begin{tabular}{@{}ccrrcr@{}}
		\toprule
		\multirow{2}{*}{Dimension}                 & Degree $p$         
            & \multicolumn{2}{c}{Quadrature points} & \multirow{2}{*}{\shortstack{Relative savings}} 
            &  \multirow{2}{*}{Restarts}\\[0.8em]
		                                          &                    & Gaussian ($q'$)       & Optimised ($q$)         \\ \midrule
		\multirow{10}{*}{2D} & 1 &   4 &   4 &  0.0\% &  1 \\
		                      & 2 &   9 &   9 &  0.0\% &  1 \\
		                      & 3 &  16 &  13 & 18.8\% &  2 \\
		                      & 4 &  25 &  19 & 24.0\% &  2 \\
		                      & 5 &  36 &  27 & 27.8\% & 98 \\
                              & 6 &  49 &  36 & 26.5\% &  5 \\
                              & 7 &  64 &  46 & 28.1\% & 12 \\
                              & 8 &  81 &  58 & 28.4\% & 23 \\
                              & 9 & 100 &  71 & 29.0\% & 3072 \\
		                      &10 & 121 &  85 & 29.8\% & 103 \\ \midrule
		\multirow{6}{*}{3D}  & 1 &   8 &   8 &  0.0\% &  1 \\
		                      & 2 &  27 &  25 &  7.4\% & 39 \\
		                      & 3 &  64 &  43 & 32.8\% & 12 \\
		                      & 4 & 125 &  74 & 40.8\% & 203 \\
		                      & 5 & 216 & 118 & 45.4\% & 517 \\
		                      & 6 & 343 & 178 & 48.1\% & 770 \\
		\bottomrule
	\end{tabular}
\end{table}


\section{Case Studies}
\label{sec:CaseStudies}

To validate the accuracy and robustness of the proposed quadrature rules, we present a set of representative case studies for solving the Poisson equation using high-order finite elements on both standard and curved geometries. We compare the performance of our optimised quadrature rules against classical Gaussian quadrature in terms of $L^2$ error convergence. In each case, we employ a quadrature rule to integrate exactly the mass matrix, and thus, its gradient when the mappings are affine. Our focus is on serendipity (trunk) spaces with varying polynomial degrees in both 2D and 3D settings.

\subsection{Test Case 1: Regular Geometries}

We begin by considering the classical Poisson problem:
\begin{equation}
    -\Delta u = f \quad \text{in } \Omega,
\end{equation}
where the domain $\Omega$ is the unit square $[0,1]^2$ in two dimensions, and the unit cube $[0,1]^3$ in three dimensions. Homogeneous Dirichlet boundary conditions are imposed throughout.

\begin{figure}[h!]
     \centering
     \begin{subfigure}[b]{0.48\textwidth}
         \centering
         \includegraphics[width=\textwidth]{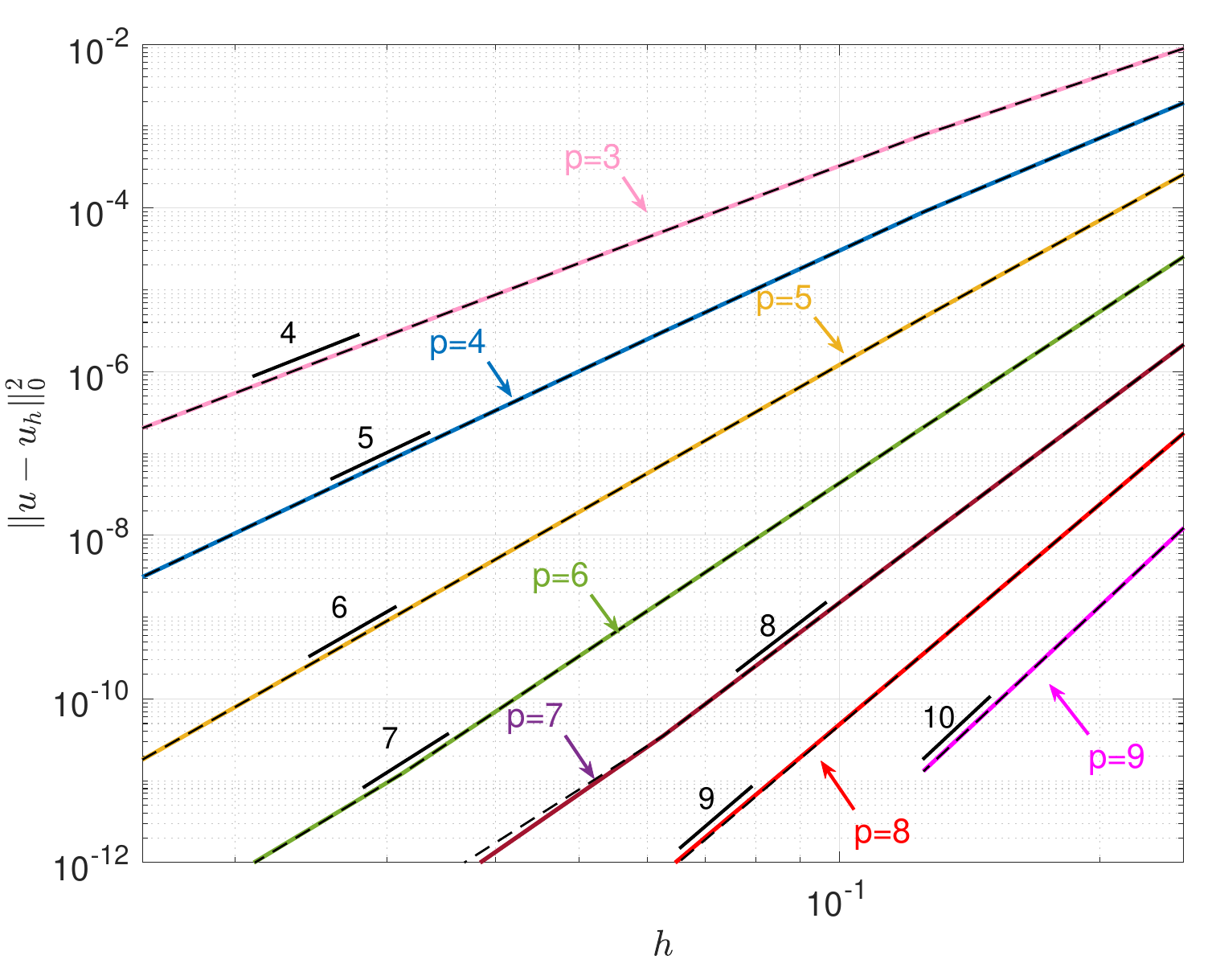}
         \caption{2D domain}
     \end{subfigure}
     \hfill
     \begin{subfigure}[b]{0.48\textwidth}
         \centering
         \includegraphics[width=\textwidth]{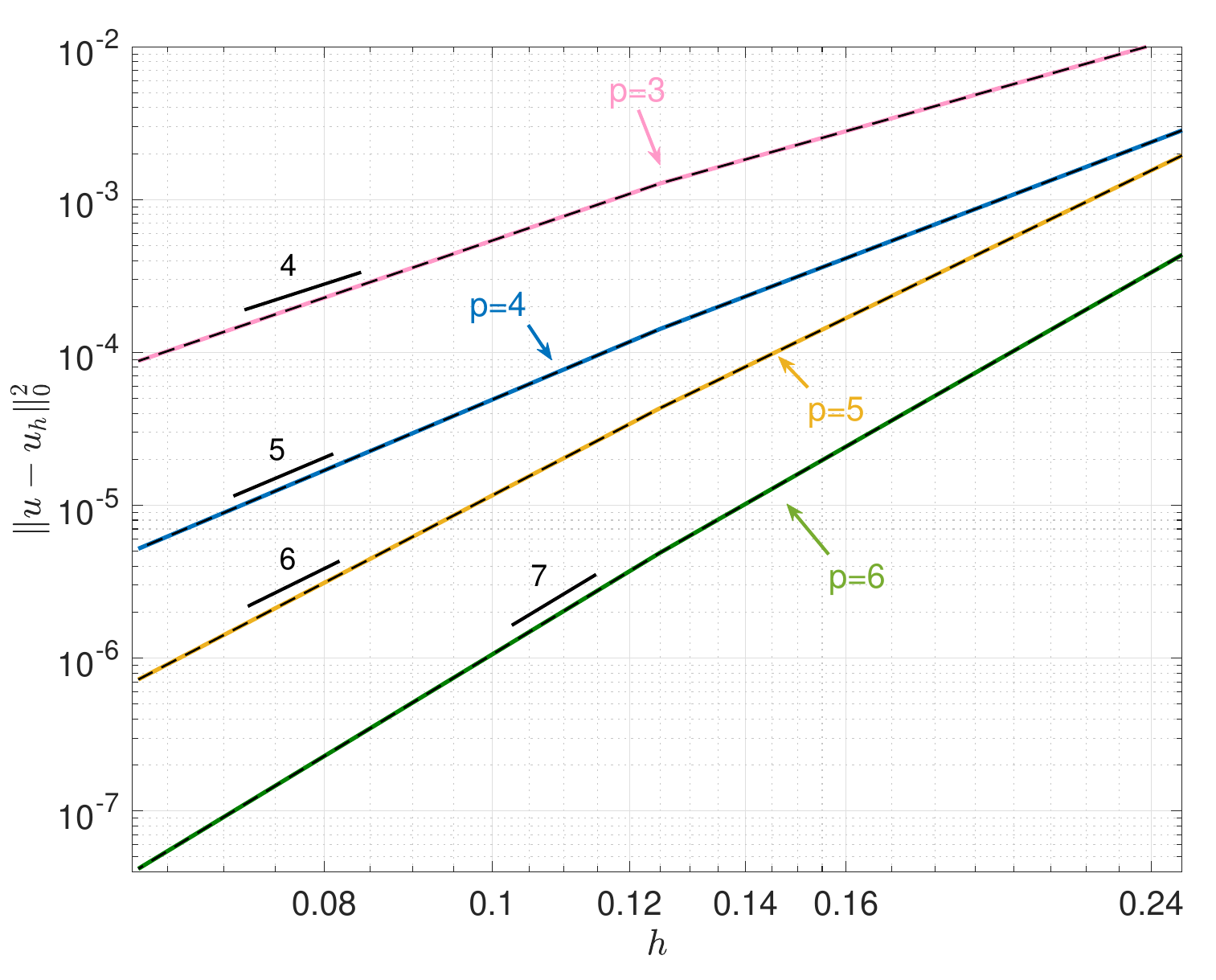}
         \caption{3D domain}
     \end{subfigure}
        \caption{$L^2$ error convergence for the proposed quadrature in 2D (left) and 3D (right). Solid lines: proposed rule; dashed lines: Gaussian rule.}
        \label{fig:conv}
\end{figure}
The weak formulation seeks a function \( u \in H^1(\Omega) \) satisfying
\begin{equation}
    (\nabla u, \nabla v) = (f, v) \quad \forall v \in H^1(\Omega).
    \label{eq:Poisson}
\end{equation}

We assess precisely numerical errors by prescribing a manufactured solution:
\begin{align}
    u(x,y) &= \sin(2\pi x)\sin(2\pi y), \quad &\text{(2D)}, \\
    u(x,y,z) &= \sin(2\pi x)\sin(2\pi y)\sin(2\pi z), \quad &\text{(3D)}.
    \label{eq:analytical}
\end{align}
The corresponding source term \( f \) and boundary data are derived accordingly.
Figure~\ref{fig:conv} illustrates the convergence in the $L^2$ norm as the mesh is refined. Solid lines denote the performance of the proposed quadrature rules, whereas dashed lines correspond to classical Gaussian integration. In both 2D and 3D cases, the proposed rules yield optimal convergence rates of order \( \mathcal{O}(h^{p+1}) \), where \( p \) is the polynomial degree. This validates their theoretical correctness and numerical consistency.
Despite using significantly fewer integration points, the proposed rules match the accuracy of classical Gaussian quadrature. This property renders them a computationally efficient alternative for high-order simulations on regular domains.

\begin{figure}[h!]
     \centering
     \begin{subfigure}[b]{0.4\textwidth}
         \centering
         \includegraphics[width=\textwidth]{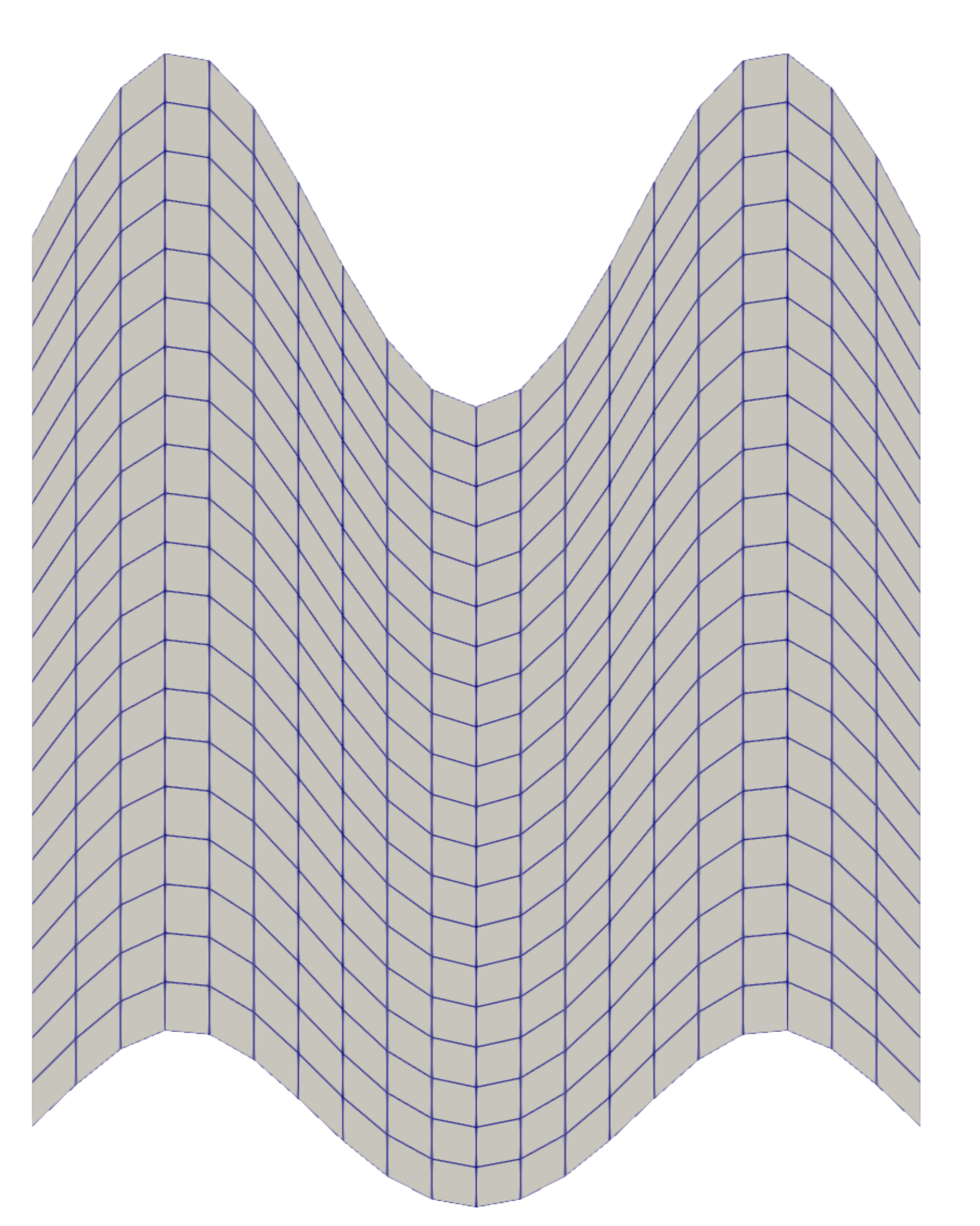}         
     \end{subfigure}
     \hfill
     \begin{subfigure}[b]{0.48\textwidth}
         \centering
         \includegraphics[width=\textwidth]{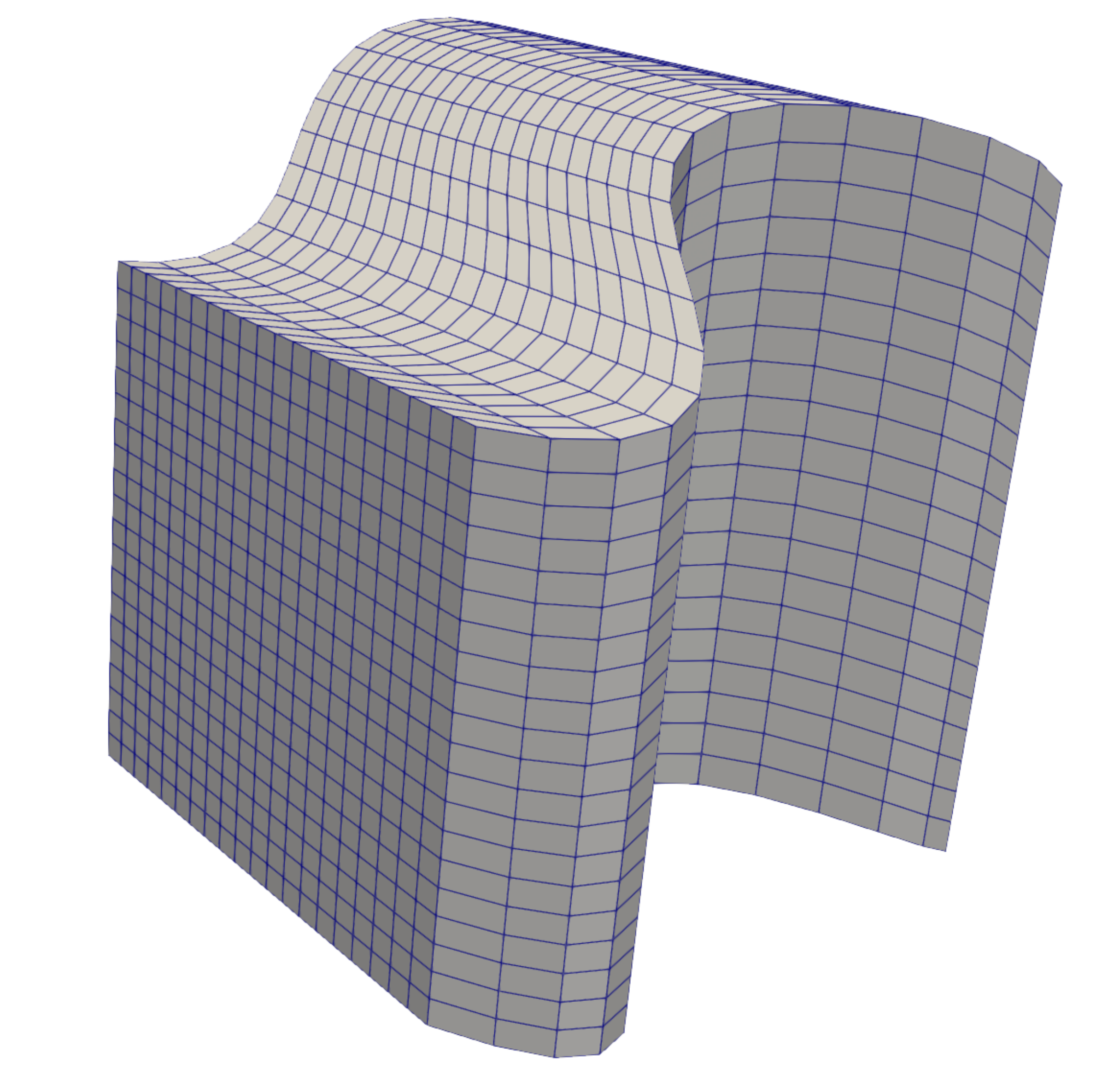}
     \end{subfigure}
        \caption{Illustration of deformed 2D (left) and 3D (right) domains.}
        \label{fig:curved_domain}
\end{figure}

\subsection{Test Case 2: Curved Geometries}

We now consider geometrically perturbed domains to evaluate the robustness of the proposed rules under non-affine transformations. The Poisson equation~\eqref{eq:Poisson} is solved on smoothly deformed versions of the unit square (2D) and unit cube (3D). In 2D, we displace the upper edge of the square vertically by \( 0.1 \sin(3\pi y) \), producing a wavy boundary and introducing an additional deformation in the $z$-direction: \( 0.1 \sin(1.5\pi y) \) for the three-dimensional case. These mappings result in curved quadrilateral and hexahedral elements. Figure~\ref{fig:curved_domain} shows representative meshes.
\begin{figure}[h!]
     \centering
     \begin{subfigure}[b]{0.48\textwidth}
         \centering
         \includegraphics[width=\textwidth]{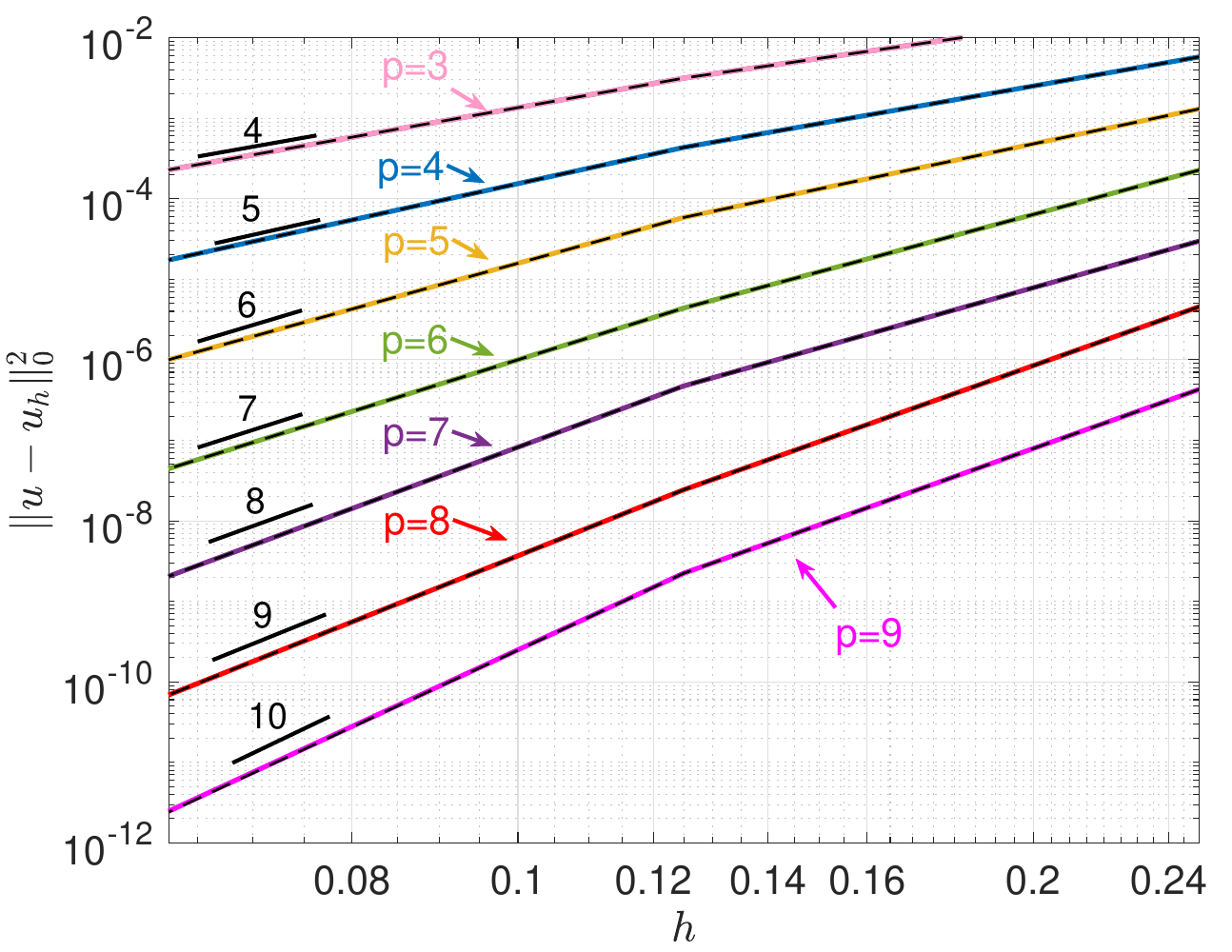}
         \caption{2D curved domain}
     \end{subfigure}
     \hfill
     \begin{subfigure}[b]{0.48\textwidth}
         \centering
         \includegraphics[width=\textwidth]{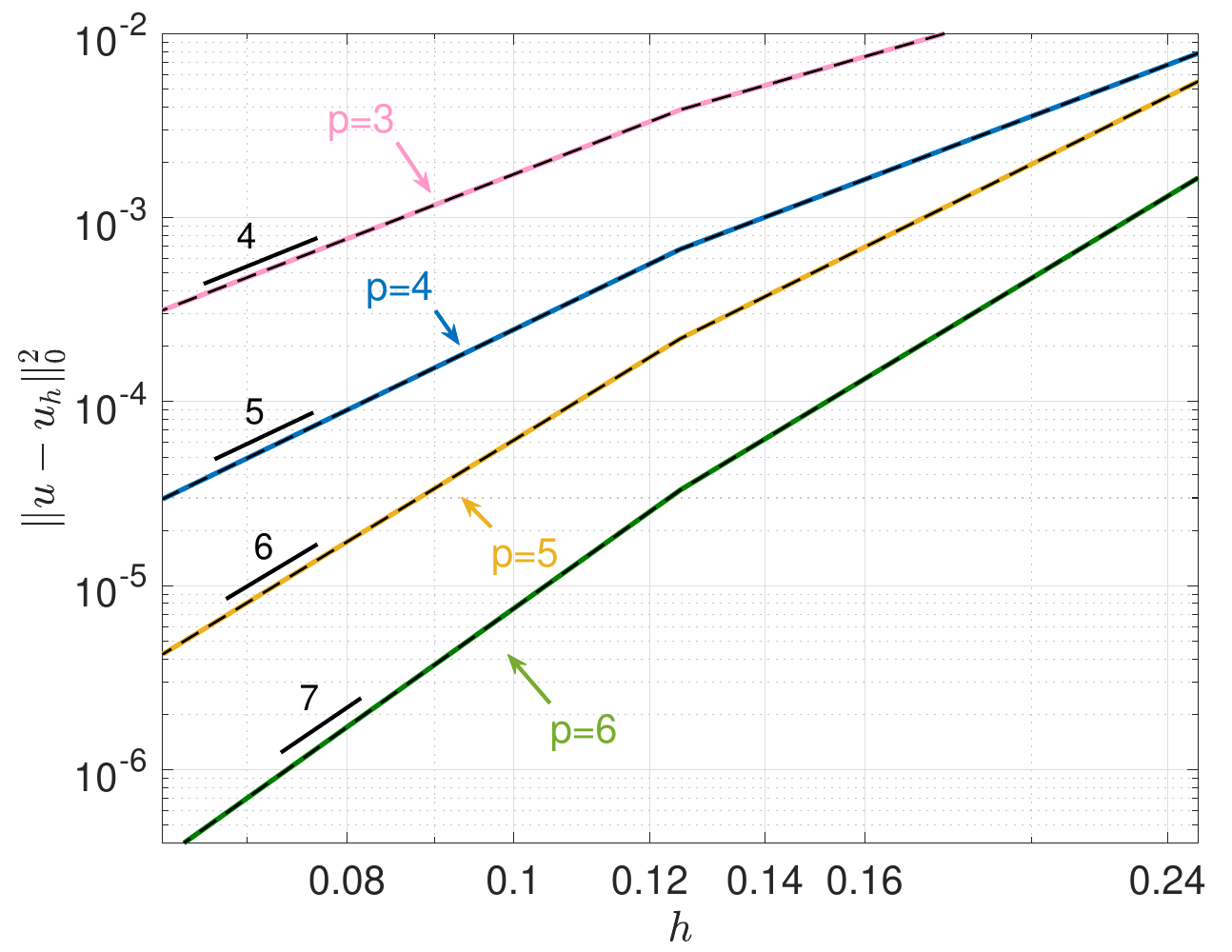}
         \caption{3D curved domain}
     \end{subfigure}
        \caption{$L^2$ error convergence on curved geometries. Solid lines: proposed quadrature; dashed lines: Gaussian quadrature.}
        \label{fig:curved_conv}
\end{figure}
Figure~\ref{fig:curved_conv} presents the convergence results for these curved configurations. Solid and dashed lines correspond to our and Gaussian quadratures, respectively. In both scenarios, the convergence remains optimal, without instabilities or loss of accuracy.

\begin{table}[h!]
    \centering
    \begin{tabular}{cccc}
        \hline
        \( p \) & \( 4^3 \) elements & \( 8^3 \) elements & \( 16^3 \) elements \\
        \hline
        2 & $3.589 \times 10^{-7}$ & $5.791 \times 10^{-6}$ & $4.107 \times 10^{-7}$ \\
        3 & $3.118 \times 10^{-8}$ & $7.199 \times 10^{-9}$ & $1.399 \times 10^{-10}$ \\
        4 & $1.852 \times 10^{-7}$ & $1.343 \times 10^{-9}$ & $5.141 \times 10^{-12}$ \\
        5 & $8.397 \times 10^{-9}$ & $4.246 \times 10^{-10}$ & $1.502 \times 10^{-12}$ \\
        6 & $2.214 \times 10^{-9}$ & $1.550 \times 10^{-11}$ & $6.376 \times 10^{-14}$ \\
        \hline
    \end{tabular}
    \caption{Absolute $L^2$ error differences between the proposed and Gaussian quadrature rules in 3D curved domains, across varying polynomial degrees and mesh sizes.}
    \label{tab:quadrature_error}
\end{table}

\medskip

The data confirms that our quadrature scheme is accurate and stable under complex geometric mappings. However, a minor discrepancy arises in the pre-asymptotic regime. The proposed method may exhibit a slightly higher error than Gaussian quadrature for low polynomial degrees or coarse meshes. This effect diminishes as the mesh is refined or the polynomial degree increases, and both methods eventually reach the same asymptotic convergence behaviour.
Table~\ref{tab:quadrature_error} quantifies this discrepancy by reporting the absolute difference in $L^2$ error for the 3D curved domain. The differences are negligible and decay rapidly with increasing resolution or polynomial degree.

Thus, these quadrature rules have excellent numerical properties across various scenarios, from regular meshes to geometrically intricate configurations. They match the convergence properties of classical Gaussian integration while requiring significantly fewer integration points, offering a promising and efficient alternative for high-order finite element computations.


\section{Conclusions}
\label{sec:Conclusions}

We address a longstanding inefficiency in the numerical integration of high-order finite element computations by proposing a systematic method to reduce the number of integration points required for exactness. In standard practice, quadrature rules for two- and three-dimensional problems are constructed as tensor products of one-dimensional Gaussian rules. While these tensor-product constructions ensure exactness for full polynomial spaces, they are often excessive for subspaces commonly used in finite element method (FEM) implementations, such as trunk or serendipity spaces. These subspaces discard many higher-order basis functions in full tensor-product spaces, rendering standard quadrature unnecessarily costly.

Following the idea behind using trunk spaces, we propose constructing minimal quadrature rules as a nonlinear optimisation problem to overcome this inefficiency. The objective is to determine the positions and weights of integration points that yield exact integration for a specified polynomial subspace. The number of equations is determined by the cardinality of the basis functions spanning the target space, while the degrees of freedom correspond to the coordinates and weights of the quadrature points. Solving this system requires a balance between accuracy and computational feasibility, especially in three-dimensional cases. Our approach employs gradient-based optimisation techniques combined with random restarts, a typical strategy in machine learning, to navigate the non-convex landscape and avoid suboptimal local minima.

Central to this strategy is trunk spaces, which significantly reduce the dimensionality of the integration problem without compromising approximation quality. By restricting the integration domain to a subset of the tensor-product basis, we lower the number of required quadrature points while ensuring exactness for the intended function space. This makes the optimisation problem tractable and yields considerably more efficient integration schemes.

Among the most compelling outcomes of our investigation is the observation that, for all practical configurations tested, the set of exact quadrature rules is potentially infinite. This multiplicity of solutions introduces a valuable degree of flexibility. It suggests one can generate families of equivalent rules, possibly optimised for different criteria such as numerical stability, locality, or computational load distribution. This insight lays the groundwork for future stochastic and randomised integration strategies developments. Such approaches could be particularly advantageous in high-performance computing environments, where load balancing and parallel efficiency are critical.

We also provide a collection of precomputed quadrature rules that are exact for trunk polynomial spaces in two and three dimensions, covering polynomial degrees up to $p=10$ in 2D and $p=6$ in 3D. These rules substantially reduce the required integration points compared to standard Gaussian tensor-product quadrature. We also introduce a flexible and general optimisation-based framework that can be used to compute exact quadrature rules for arbitrary polynomial spaces, including those with anisotropy or non-uniform degree distributions. By synthesising mathematical understanding of polynomial structure with tools from numerical optimisation, we have created a practical and extensible methodology for improving the efficiency of numerical integration in high-order FEM contexts. These advances are not merely technical; they represent a step forward in enabling more adaptive, scalable, and responsive computational models. Thus, we improve the performance of numerical integration schemes and offer a conceptual shift toward more adaptable and intelligent quadrature design in scientific computing.

In the future, we will explore extending this method to accommodate dynamically changing geometries, adaptive meshing strategies, and real-time simulation requirements. We also foresee applications in stochastic FEM, where the quadrature rules' non-uniqueness may be exploited to inject beneficial randomness into the integration process, enhancing robustness and parallelism.


\section*{Declaration of generative AI and AI-assisted technologies in the writing process}

During the preparation of this work the author(s) used ChatGPT in order to check the grammar style and wording of specific paragraphs in the introduction and conclusions sections. After using this tool, the author(s) reviewed and edited the content as needed and take(s) full responsibility for the content of the publication.

\section*{Acknowledgements}

The authors would like to dedicate this work to the memory of our colleague Ali Hashemian.

T. Teijeiro is supported by the grant RYC2021-032853-I/MCIN/AEI/10.13039/501100011033 funded by the Spanish Ministry of Science and Innovation and by the European Union NextGenerationEU/PRTR.

David Pardo has received funding from the following Research Projects/Grants: 
European Union's Horizon Europe research and innovation programme under the Marie Sklodowska-Curie Action MSCA-DN-101119556 (IN-DEEP). 
PID2023-146678OB-I00 funded by MICIU/AEI /10.13039/501100011033 and by FEDER, EU;
"BCAM Severo Ochoa" accreditation of excellence CEX2021-001142-S funded by MICIU/AEI/10.13039/501100011033; 
Basque Government through the BERC 2022-2025 program;
Consolidated Research Group MATHMODE (IT1456-22) given by the Department of Education of the Basque Government; 
BCAM-IKUR-UPV/EHU, funded by the Basque Government IKUR Strategy and by the European Union NextGenerationEU/PRTR;


\bibliography{bibliography}

\end{document}